\documentclass{article}

\usepackage[english,russian]{babel}

\usepackage{amsmath}
\usepackage{amsfonts}
\usepackage{amssymb,latexsym}
\usepackage{amsthm}

\newtheorem{theorem}{Теорема}[section]
\newtheorem{lemma}{Лемма} [section]
\newtheorem{remark}{Замечание}
\newtheorem{corollary}{Следствие}

\newtheorem{{thebibliography}}{REf}

\begin{document}

\author{М.З. Двейрин}

\title{О скорости полиномиальной аппроксимации целых функций и их свойствах}



\maketitle

Аннотация. В статье выясняется связь между порядком и типом целой
функции и скоростью наилучшей полиномиальной аппроксимации для большого семейства
банаховых пространств функций, аналитических в единичном  круге. Найдены соотношения, определяющие порядок и тип
целой функции через последовательность  ее наилучших приближений.
Полученные результаты являются обобщением более ранних результатов
Редди, И.И. Ибрагимова и Н.И. Шихалиева, С.Б. Вакарчука, Р. Мамадова. В тех же банаховых пространствах аналитических функций рассмотрены также некоторые вопросы приближения полиномами с целыми коэффициентами.

Annotation. The paper explores connection between the order and the type of an entire function and the speed of the best polynomial approximation in the unit disk. The relations which define the order and the type of an entire function through the sequence of its best approximations, have been found. The results were obtained by generalization previous results of Reddy, I.I. Ibragimov and N. I. Shyhaliev, S. B. Vakarchyk, R. Mamadov.

2000 MSC.  41А10, 41А25, 41А58.

Ключевые слова и фразы. Целая функция, наилучшее приближение,
порядок целой функции, тип целой функции.

\section{Введение}

Рассмотрим банахово пространство $X$, образованное
аналитическими в единичном круге  $\mathbb{D}$ функциями, имеющими конечную норму $\|\cdot\| $\,. \, Будем
считать, что $\| \cdot \| $\, помимо обычных свойств
нормы удовлетворяет также условиям

\begin{equation}
\label{trivial} i) \quad \, \|f(\cdot e^{it})\| \, \equiv
\|f(\cdot)\| \,
\end{equation}
для всех $t \in \mathbb{R}$ и $f \in X;$

\begin{equation}
ii) \quad \, \|f(\cdot)\| \,<\infty  \,
\end{equation}
для   любой целой функции (т.е. пространство $X$ содержит все целые функции);

\begin{equation}
iii) \quad \, \|  \frac {1}{2\pi}\int\limits_{0}^{2\pi}f(ze^{it})g(t)
\,dt \| \leq \, \frac {1}{2\pi}\int\limits_{0}^{2\pi}|g(t)|
\,dt \,\, \|f(\cdot)\|
\end{equation}
для   любых функций $f \in X$  и  $g \in L[0; 2\pi]$   (иначе говоря, для любых $f \in X$  и  $g \in L[0; 2\pi]$ \, $\|f \ast g \| \leq \| f \|\, \|g\|_{L[0; 2\pi]}$).

Этим требованиям удовлетворяет норма в целом ряде функциональных
пространств, являющихся объектом многочисленных исследований (автору неизвестны примеры пространств, в которых выполняются условия i),  ii) и при этом не выполняется условие iii) ).
Приведем некоторые из них.

 1) Пространство $B$ функций, аналитических в единичном круге $\mathbb{D}$ и
непрерывных на его замыкании $\overline {\mathbb{D} }$ с нормой

\[
\|f\|=\max\limits_{z\in
\overline{\mathbb{D}}}|f(z)|<\infty \,.
\]

 2) Пространства Харди $H_p$ ($p \geq 1$) функций, аналитических в круге
$\mathbb{D}$   с нормой

\[
\|f\|=\sup \limits_{0<r<1}M_{p} \,(f,r),\quad M_{p}
\,(f,r):= \left(\frac {1}{2\pi}\int\limits_{0}^{2\pi}|f(re^{it})|^p
\,dt\right)^\frac{1}{p}\;,\quad p\in [1; \infty);
\]

\[
\|f\|=\sup\limits_{z\in \mathbb{D}}|f(z)|\,,\, \quad
p=\infty.
\]

 3) Пространства Бергмана $H_p^\prime$ функций, аналитических в круге
$\mathbb{D}$ при $p\in [1; \infty)$ с нормой

\[
\|f\|=\left(\frac {1}{\pi} \int \limits_{\;\; z \in
D}\int |f(x+iy)|^p \;dxdy\right)^\frac{1}{p}\;,
\]
и обобщенные (весовые) пространства Бергмана $H_{p, \, \rho}^\prime$  функций, аналитических в круге
$\mathbb{D}$ при $p\in [1; \infty)$ с нормой
\[
\|f\|=\left(\frac {1}{\pi} \int \limits_{\;\; z \in
D}\int |f(x+iy)|^p \;\rho(|z|)dxdy\right)^\frac{1}{p}\;
\]
и радиальным весом $\rho(|z|)$.

 4) Пространства  $A_p, \quad p \in (0;1)\; $ функций,
аналитических в круге $\mathbb{D}$  с нормой

\[
\|f\|=\int\limits_{0}^{1}(1-r)^{\frac{1}{p}-2}M_{1} \,(f,r)
\,dr\, ,
\]
впервые изучавшиеся Харди и Литллвудом \cite{HL2} и позднее Ромбергом, Дюреном и Шилдсом  \cite{DRS}.

 5) Пространства  $\mathcal{B}_{p,\,q,\,\lambda}\, , \quad
0<p<q\leq\infty, \quad \lambda > 0,\; $ функций,
аналитических в круге $\mathbb{D}$ с нормой

\[
\|f\|=\left\{ \int \limits_{0}^{1}(1-r)^{\lambda
 \, p \;q \,(q-p)^{-1}}M_q^\lambda(f,r) \,dr\right\}^\frac{1}{\lambda},\, \quad
\lambda<\infty,
\]

\[
\|f\|=\sup\limits_{0<r<1}\left\{ (1-r)^{
 \, p \;q \,(q-p)^{-1}}M_{q}\,(f,r)\right\} \, ,\quad \lambda=\infty,
\]
введенные {Харди и Литллвудом в работе \cite{HL2} (см. также \cite{Gvar}).

 6) Пространства со смешанной нормой $H^{p,\,q,\,\alpha} , \,
(p, q \geq 1, \quad\,  \alpha >0),$ образованные функциями,
аналитическими в круге $\mathbb{D}$ с конечной нормой

\[
\|f\|=\left\{ \int \limits_{0}^{1}(1-r)^{q\alpha -1}M_p^q(f,r) \,dr\right\}^\frac{1}{q} \,
,    \quad q<\infty,
\]

\[
\|f\|=\sup\limits_{0<r<1}\left\{ (1-r)^{
 \alpha}M_{p}\,(f,r)\right\} \, ,\quad q=\infty,
\]
введенные Харди и Литллвудом в работе  \cite{HL2}.    Заметим, что пространства   $H^{p,\,q,\,\alpha}$  и
$\mathcal{B}_{p,\,q,\,\lambda}\,$  отличаются лишь способом введения параметров.

 7) Пространство $\quad BMOA$ \cite{Swed}, состоящее из функций $f\in H_{1}$ с нормой
\[
\|f\|=\sup\limits_{I}\int \limits_{I} |f(\zeta)-f_{I}|d\sigma(\zeta)\,,
\]
где $f(\zeta)$ - граничные значения функции $f(z)$ на единичной окружности,
а $f_{I}$ - среднее арифметическое значение функции $f(\zeta)$ на дуге $I$.

 8) Пространства типа Блоха $\mathcal{B}_{\alpha}\,$,  $\alpha\in(0, \infty)$, состоящие из функций, аналитических в  $D$ с конечной нормой
\[
\|f\|=|f(0)|+ \sup\limits_{z\in D} (1-|z|^{2})^{\alpha} |f^{\prime}(z)|\,.
\]
Пространства  $\mathcal{B}_{\alpha}\,$  являются банаховыми  \cite{Zhu},  при $\alpha=1$   $\mathcal{B}_{\alpha}$ совпадает с пространством Блоха $\mathcal{B}$.

 9) Введенные Е.М. Дынькиным \cite{Dink} пространства  $\mathcal{A}_{p,q}^s (\mathbb{D})\,$
 функций, являющиеся аналогами классов О.В. Бесова $\mathcal{B}_{p,q}^s
[-1;\;1]\,$.  Эти пространства образованы функциями $f\in H_{p}$, $p\in [1; \infty]$  с нормой

\[
\|f\|=\left\{ \int \limits_{0}^{1}\left (  \frac{\omega_{m}(f, \, t)_{p}}{t^{s}} \right)^{q}\frac{dt}{t}\right\} ^\frac{1}{q} \, +  \sup \limits_{0<r<1}M_{p} \,(f,r).
\]

Здесь  $q\in [1; \infty], \, s>0, \, m>s$ - натуральное число,  $\omega_{m}(f, \, t)_{p}$ - $m$-ый модуль гладкости в пространстве  $L_p\,$ функции  $f(e^{i\cdot})$, представляющей собой радиальные предельные значения $f$.  Случай $q=\infty$ трактуется традиционно.

10) Обобщенные пространства  Дирихле $\mathcal{D}_{p}(\alpha)\,$ функций,
аналитических в  $\mathbb{D}$,  с нормой
\[
\|f(z)\|=\left( \sum\limits_{k=0}^{\infty}|c_{k}|^{p}\,\alpha_{k}\right)^{1/p}\,,
\]
где   $ c_{k}=c_{k}(f)$ - коэффициенты Тейлора функции $f$,  $p\geq 1$,  ${\alpha}=\{\alpha_{k}\}$ - фиксированная последовательность положительных чисел с условиями
$$ \limsup \limits_{k\rightarrow\infty}\left(\alpha_{k}\right)^{\frac{1}{k}}< \infty, \quad
 \liminf \limits_{k\rightarrow\infty}\left(\alpha_{k}\right)^{\frac{1}{k}}\geq 1 . $$

Отметим, что приведенные выше примеры функциональных  пространств со свойствами i), ii) и iii) не исчерпывают их многообразия.

Обозначим  $E_n(f) \equiv E_n(f ,\,L_n)\;$
наилучшее приближение функции  $f\in X$ элементами линейного
подпространства $L_n$:

\[
E_n(f):= \inf \limits_{p\in L_n}\|f-p \, \|\;.
\]

\noindent В качестве аппроксимирующего подпространства $L_n$ мы
будем  рассматривать совокупность $\mathcal{P}_n$ алгебраических
полиномов комплексной переменной степени не выше $(n-1)$,
$p^{*}_{n}$ - полином наилучшего приближения степени не выше  $(n-1)$ для функции $f$,
т.е. такой, что  $E_n(f)= \|f-p^{*}_{n}\|$. В статье
\cite{Rd1} найдены соотношения, определяющие порядок и тип целой
функции через последовательность $E_n(f)$ ее наилучших
приближений в случае $X=H_2^\prime$. В  1976г. Ибрагимов И.И. и
Шихалиев Н.И. получили подобные соотношения для пространств
$X=H_p^\prime$ с произвольным $p\geq 1$. Приведем формулировки
теорем из их работы \cite{Ibr1} (см. также \cite{Ibr2}).

\begin{theorem}\label{sol1}
      Для того, чтобы    $f(z)\in H_p^\prime$, ($p\geq1$, любое)
      была целой функцией, необходимо и достаточно, чтобы

      \begin{equation}
      \lim_{n{\rightarrow} \infty}
        (E_n(f))^{\frac{1}{n}}=0.
      \label{11}
      \end{equation}
      \end{theorem}

\begin{theorem}\label{sol2}
      Для того, чтобы    $f(z)\in H_p^\prime$, ($p\geq1$, любое)
      была целой функцией конечного порядка $\rho$, необходимо и достаточно, чтобы

      \begin{equation}
        \limsup \limits_{n\rightarrow\infty}
        \,\frac{n \ln n}{-\ln E_n(f)}=\rho.
      \label{12}
      \end{equation}
      \end{theorem}

\begin{theorem}\label{sol3}
      Для того, чтобы    $f(z)\in H_p^\prime$, ($p\geq1$, любое)
      была целой функцией конечного порядка $\rho$ и нормального типа $\sigma$, необходимо и достаточно, чтобы
      \begin{equation}
       \limsup \limits_{n\rightarrow\infty}
        \,n( E_n(f))^{\frac{\rho}{n}}=\sigma e\rho.
      \label{13}
      \end{equation}
      \end{theorem}

   В 1990г. в \cite{Wac1, Wac2} аналогичные результаты получены С.Б.
Вакарчуком для пространств $\mathcal{B}_{p,\,q,\,\lambda}$.
Р. Мамадов в своей кандидатской диссертации (Душанбе, 2009г., \cite{Mamad})
распространил эти результаты на весовые пространства Бергмана с радиальным весом $\rho(|z|)$.
В настоящей статье эти результаты распространяются на широкую
совокупность нормированных пространств функций, аналитических в
единичном круге, включающую среди прочих приведенные выше
пространства 1) $-$ 10).  Метод доказательства фрагментами повторяет  доказательства в указанных выше
 работах, но в целом не совпадает ни с одним из них.

Помимо этого в  статье рассмотрены некоторые вопросы
приближения функций в пространстве $X$ полиномами $p_{n}(z)$ с
целыми комплексными коэффициентами (считаем коэффициенты целыми,
если у них действительная и мнимая части являются целыми числами;
совокупность комплексных полиномов степени не выше $(n-1)$ с целыми
коэффициентами   будем обозначать  ${\mathcal{P}}_{n}[\mathbb{Z}]$ ).
Полученные результаты приведены в разделе 5. С результатами
предыдущих разделов их объединяет метод доказательства, основанный
на одних и тех же леммах.

\section{ Формулировка результатов.}

Здесь  будут приведены формулировки основных результатов статьи.

\begin{theorem}\label{sol21}
Пусть  $f\in X$. Условие
 \begin{equation}
\lim\limits_{n\rightarrow
\infty}(E_{n}(f))^{\frac{1}{n}}=0
 \end{equation}
является необходимым и достаточным для того, чтобы функция $f$ была
целой.
           \label{21}
      \end{theorem}

\begin{theorem}\label{sol22}
Для того, чтобы  функция $f\in X$ была целой конечного порядка
$\rho\in(0;\infty)$ необходимо и достаточно, чтобы
\begin{equation}
\label{L11}\limsup \limits_{n\rightarrow\infty}\,\frac{n\ln
n}{\ln\frac{\|z^{n}\|}{E_{n}(f)}}=\rho .
\end{equation}
 \label{22}
    \end{theorem}

\begin{theorem}\label{sol23}
Пусть существует конечный предел $\lim\limits_{n\rightarrow \propto}
(\|z^{n}\|)^\frac{1}{n}=\mu $. Для того, чтобы  функция $f\in
X$ была целой функцией конечного порядка $\rho\in(0;\infty)$ и
нормального типа $\sigma\in(0;\infty)$ необходимо и достаточно,
чтобы
\begin{equation}
\label{L15}\limsup\limits_{n\rightarrow\infty}
\,\frac{n}{e\rho}\left(\frac{E_{n}(f)}{\|z^{n}\|}\right)
^{\frac{\rho}{n}}=\sigma .
\end{equation}
\label{23}
 \end{theorem}

\section{ Применение к конкретным пространствам.}
Выполнение условий i) и ii) во всех вышеприведенных примерах пространств очевидно. Покажем, что в приведенных выше функциональных пространствах выполняется  также условие iii).

\begin{lemma}
Условие iii) выполнено в пространствах   $B$;\, $H_p$ ($p \geq 1$); \, $H_{p, \, \rho}^\prime$  ($p \geq 1$); \,
$A_p, \quad (p \in (0;1))$;  \, $H^{p,\,q,\,\alpha}\, , \quad
(p,q \geq 1, \quad \alpha > 0)$; \, $\quad BMOA$; \,$\mathcal{B}_{\alpha}, \quad (\alpha \in (0;\infty))$;  $\mathcal{A}_{p,q}^s (\mathbb{D})\, \,(p,q\in [1; \infty], \,s>0) $; \,  $\mathcal{D}_{p}(\alpha)\,$ ($p \geq 1$).

 \begin{proof}
В пространствах $B$ и $A_{p}$  выполнение свойства  iii) очевидно;  в пространствах  $H_p$, $H_{p, \, \rho}^\prime$  при $p \geq 1$ свойство  iii) очевидным образом следует из обобщенного неравенства Минковского для пространств $L_{p}$.
В случае пространства $\quad BMOA$

\[
\|  \frac {1}{2\pi}\int\limits_{0}^{2\pi}f(ze^{it})g(t)\,dt \|=
\]

\[
 = \sup\limits_{I}\int \limits_{I} \left| \frac {1}{2\pi}\int\limits_{0}^{2\pi}f(e^{i(t+\varphi)})g(t)dt -
  \frac {1}{|I|}\int\limits_{I} \left( \frac {1}{2\pi}\int\limits_{0}^{2\pi}f(e^{i(t+u)})g(t)dt \right) du \right| d\varphi=
\]
\[
= \sup\limits_{I}\int \limits_{I} \left| \frac {1}{2\pi}\int\limits_{0}^{2\pi}f(e^{i(t+\varphi)})g(t)
\,dt -  \frac {1}{2\pi}\int\limits_{0}^{2\pi} \left(\frac {1}{|I|}\int\limits_{I}f(e^{i(t+u)})g(t)du\right)dt \right|\, d\varphi=
\]
\[
= \sup\limits_{I}\int \limits_{I} \left| \frac {1}{2\pi}\int\limits_{0}^{2\pi}g(t)\left(f(e^{i(t+\varphi)})
 -  \frac {1}{|I|}\int\limits_{I}f(e^{i(t+u)})du\right)dt \right|\, d\varphi \leq
 \frac {1}{2\pi}\int\limits_{0}^{2\pi}|g(t)|
\,dt \,\, \|f(\cdot)\|
\]
ввиду свойства i) нормы в пространстве $X.$

Убедимся в выполнении условия  iii) в пространстве  $\mathcal{B}_{\alpha}$.
\[
\|f\ast g\|= \left|\frac {1}{2\pi}\int\limits_{0}^{2\pi} f(0)g(t)\,dt\right| + \sup\limits_{z\in D} (1-|z|^{2})^{\alpha}\, \left|\frac {1}{2\pi}\int\limits_{0}^{2\pi} f^{\prime}(ze^{it})g(t)e^{it}\,dt \right|\,\leq
\]

\[
\leq |f(0)| \, \|g(t)\|_{L}\, + \sup\limits_{z\in D}\,  \frac {1}{2\pi}\int\limits_{0}^{2\pi} (1-|z|^{2})^{\alpha}\, |f^{\prime}(ze^{it})g(t)e^{it}|\,dt \,\leq
\]

\[
\leq \|g(t)\|_{L} (|f(0)| + \sup\limits_{z\in D}\, (1-|z|^{2})^{\alpha}\, |f^{\prime}(ze^{it})|)= \|f\|\, \|g(t)\|_{L} .
\]

Пусть теперь  $f \in H^{p,\,q,\,\alpha}$,  $g \in L_{[0;\;2\pi]}, \quad z=re^{i\varphi}$.  В пространстве  $H^{p,\,q,\,\alpha}$ используя обобщенное неравенство Минковского имеем

\[
\|f\ast g\|=\left\{ \int \limits_{0}^{1}(1-r)^{q\alpha -1} \left( \frac {1}{2\pi}\int\limits_{0}^{2\pi} | \frac {1}{2\pi}\int\limits_{0}^{2\pi} f(re^{i(\varphi+t)})g(t)\,dt|^{p} d\varphi \right)^{\frac{q}{p}} \right\}^{\frac{1}{q}}\leq
\]

\[
\leq \left\{ \int \limits_{0}^{1}(1-r)^{q\alpha -1} \left( \frac {1}{2\pi}\int\limits_{0}^{2\pi} \,dt \left( \frac {1}{2\pi}\int\limits_{0}^{2\pi} |f(re^{i(\varphi+t)})g(t)|^{p} d\varphi \right)^{\frac{1}{p}} \right)^{q} \right \}^{\frac{1}{q}} = \|f\|\, \|g\|_{L}.
\]

Докажем, что  свойство  iii) выполнено и в пространстве  $\mathcal{A}_{p,q}^s (\mathbb{D})\,$.

Из интегрального представления свертки
\[
f\ast g \,\,\,(e^{i\varphi})= \frac {1}{2\pi}\int\limits_{0}^{2\pi} f(e^{i(\varphi+t)})g(t)\,dt,
\]
вследствие линейности оператора $m$-ой разности $\Delta^{h}_{m}(\cdot, \, \varphi)$
\[
\Delta^{h}_{m}(f\ast g, \varphi)= \frac {1}{2\pi}\int\limits_{0}^{2\pi} \Delta^{h}_{m}(f(e^{i(\cdot \,\, +t)}), \varphi) \, g(t)\,dt
\]
и, используя обобщенное неравенство Минковского, получаем

\[
\omega_{m}(f \ast g, \, h)_{p} \leq  \omega_{m}(f, \, h)_{p} \, \|g\|_{L} .
\]
 Из этого неравенства и отмеченной раньше справедливости свойства  iii) в пространстве $H_{p} \,$  следует справедливость свойства iii) в пространстве  $\mathcal{A}_{p,q}^s (\mathbb{D})\,$.

Докажем выполнение свойства iii) в пространстве $\mathcal{D}_{p}(\alpha)\,$. Пусть  $ c_{k}$ - коэффициенты Тейлора функции $f \in \mathcal{D}_{p}(\alpha)$,  $ b_{k}$ - коэффициенты Фурье суммируемой функции $g$, $z \in \mathbb{D}$. Тогда
\[
f \ast g(z)= \frac {1}{2\pi}\int\limits_{0}^{2\pi}f(ze^{it})g(t)
\,dt = \sum\limits_{k=0}^{\infty} c_{k}b_{-k}z^{k}
\]
и
\[
\|f\ast g\|=\left( \sum\limits_{k=0}^{\infty}|c_{k}b_{-k}|^{p}\,\alpha_{k}\right)^{1/p}\,\leq \sup\limits_{k} |b_{k}| \left( \sum\limits_{k=0}^{\infty}|c_{k}|^{p}\,\alpha_{k}\right)^{1/p} \leq \|f\|\,\|g\|_{L}.
\]

\end{proof}

\end{lemma}

Убедимся, что в приведенных выше пространствах выполнено условие теоремы 2.3 на поведение последовательности $\|z^{n}\|$. Для этого покажем, что в этих пространствах  $\lim\limits_{n\rightarrow
\infty}(\|z^{n}\|)^\frac{1}{n}=1$.
В пространствах $B$ и $H_p$, ($p\geq1$) при $n\geq0 \quad$
$\|z^{n}\|=1$ и, следовательно, условия теоремы
\ref{sol23}  в этих пространствах выполнены. В
пространствах $H_p^\prime\,$, ($p\geq1$)   при $n\geq0\quad$
$\|z^{n}\|=(np+2)^{\frac{-1}{p}}$ и условия теоремы
\ref{sol23}  в этих пространствах также выполняются,
формулировки теорем \ref{sol21} - \ref{sol23} совпадают с соответствующими теоремами из
\cite{Ibr1}.  В пространствах $\mathcal{B}_p$ при $\, \, p \in (0;1)\;$
$\|z^{n}\|=(2\pi
B(\frac{1}{p}-1,\,np+1))^{\frac{1}{p}}$, где $B(\cdot,\,\cdot)$ -
бета-функция Эйлера.  Из свойств бета-функции следует, что в этом
пространстве $\lim\limits_{n\rightarrow
\infty}(\|z^{n}\|)^\frac{1}{n}=1$ и к пространствам
$\mathcal{B}_p$ теорема \ref{sol23}  применима. В
пространствах  $\mathcal{B}_{p,\,q,\,\lambda}\,$ при
$\lambda<\infty\quad$  $\|z^{n}\|= B(\lambda
n+1,\,\frac{\lambda pq}{q-p}+1)$  и, как и в предыдущем случае,
$\lim\limits_{n\rightarrow
\infty}(\|z^{n}\|)^\frac{1}{n}=1$.  При
$\lambda=\infty\quad$  $\|z^{n}\|=
\sup\limits_{0<r<1}\,r^{n}(1-r)^{\frac{ pq}{q-p}}$. Поскольку
$(1-\frac{1}{n})^{n}n^{\frac{ pq}{p-q}} \leq\|z^{n}\|<1$, то
$\lim\limits_{n\rightarrow\infty}(\|z^{n}\|)^\frac{1}{n}=1$
 и теорема \ref{sol23}  применима к пространствам
 $\mathcal{B}_{p,\,q,\,\lambda}$. В этом случае утверждения теорем 2.1 - 2.3 совпадают с
 результатами С.Б. Вакарчука \cite{Wac1}.
 В весовых пространствах Бергмана  $H_{p,\rho}^\prime$

 \[
\|z^{n}\|=\left(2 \int \limits_{0}^{1}\; t^{pn+1} \;\rho(t)dt\right)^\frac{1}{p}\;
\leq \left(2 \int \limits_{0}^{1}\; t \;\rho(t)dt\right)^\frac{1}{p}\;
\]
и  $\limsup
\limits_{n\rightarrow\infty}\,(\|z^{n}\|)^{\frac{1}{n}}\leq 1$
в случае, если функция $ t\rho(t)$ суммируема на отрезке [0; 1].
С другой стороны, при $\varepsilon\in (0; 1)$
\[
\|z^{n}\|\geq\left(2 \int \limits_{1-\varepsilon}^{1}\; t^{pn+1} \;\rho(t)dt\right)^\frac{1}{p}\;
\geq (1 - \varepsilon)^n \; \left(2 \int \limits_{1-\varepsilon}^{1}\; t \;\rho(t)dt\right)^\frac{1}{p}\;
\]
и  $\liminf
\limits_{n\rightarrow\infty}\,(\|z^{n}\|)^{\frac{1}{n}}\geq (1 - \varepsilon)$,
если вес $\rho(t)$ удовлетворяет естественному условию, что $ \int \limits_{1-\varepsilon}^{1}\; t \;\rho(t)dt > 0$
при каждом $\varepsilon\in (0; 1)$.  Ввиду произвольности $\varepsilon\in (0; 1)$ с учетом оценки $\limsup
\limits_{n\rightarrow\infty}\,(\|z^{n}\|)^{\frac{1}{n}}\leq 1$  имеем
$\lim\limits_{n\rightarrow\infty}(\|z^{n}\|)^\frac{1}{n}=1$.  Таким образом,
теорема 2.3 применима и к пространствам $H_{p,\rho}^\prime$. В этом случае
мы получаем соответствующие результаты Р. Мамадова (см. \cite{Mamad}).

Проверим,  что условия теоремы 2.3   выполнены также в пространстве  $\mathcal{A}_{p,q}^s (\mathbb{D})\,$. Поскольку  $|\Delta^{h}_{m}(e^{i n \cdot})|=|(1 - e^{inh})^{m}|=(2 \sin\frac{nh}{2} )^{m}$,  то

\[
\| z^{n}\|=\left\{ \int \limits_{0}^{\frac{\pi}{n}}\left (  \frac{(2 \sin\frac{nt}{2} )^{m}}{t^{s}} \right)^{q}\frac{dt}{t} +  \int \limits_{\frac{\pi}{n}}^{1} 2^{mq}t^{-sq-1} dt \right \} ^\frac{1}{q} \, +  1 \leq
\]

\[
\leq\left\{ \int \limits_{0}^{\frac{\pi}{n}} n^{mq} t^{(m-s)q-1}dt +  \int \limits_{\frac{\pi}{n}}^{1} 2^{mq}t^{-sq-1} dt  \right \} ^\frac{1}{q} \, +  1 \leq  C\, n^{q},
\]
где C - положительное, зависящее только от $m, s, q$. Из этой оценки следует, что  в пространстве
$\mathcal{A}_{p,q}^s (\mathbb{D})\,$ $\lim\limits_{n\rightarrow\infty}(\|z^{n}\|)^\frac{1}{n}=1$.

Перейдем теперь к пространству $\quad BMOA$.   Положим $f=z^{n}$,  найдем $f_{I}$  и оценим $\| z^{n}\|$ в $\quad BMOA$.   Пусть $I$ дуга единичной окружности с концами в точках $e^{it_{1}}$ и $e^{it_{2}}, \quad t_{2}>t_{1}, \quad x=\frac{t_{1}+t_{2}}{2}, \quad h=\frac{t_{2}-t_{1}}{2}$. Тогда
\[
f_{I}=\frac{1}{t_{2}-t_{1}} \, \int \limits_{t_{1}}^{t_{2}} e^{int} dt = \frac{\sin nh}{nh}\, e^{inx}
\]
Положим
\[
A=\frac{1}{t_{2}-t_{1}} \, \int \limits_{t_{1}}^{t_{2}} |e^{int}- \frac{\sin nh}{nh}\, e^{inx}| \,dt = \frac{1}{2h}\, \int \limits_{-h}^{h} |e^{int}- \frac{\sin nh}{nh}\, | \,dt.
\]
Поскольку $A\leq 2$, то   $\|z^{n}\|= \sup\limits_{h\in [0; \, \pi]} A \leq 2$. Оценим теперь $\|z^{n}\|$ снизу.
\[
\|z^{n}\|\geq \sup\limits_{h=\frac{\pi}{2n}} A= \sup\limits_{n} \frac{n}{\pi}
\, \int \limits_{\frac{-\pi}{2n}}^{\frac{\pi}{2n}} \left(1-  \frac{2}{\pi}\,\cos nt + \frac{4}{\pi^{2}} \right )^{\frac{1}{2}} \,dt \geq \sup\limits_{n} \frac{2n}{\pi}
\, \int \limits_{0}^{\frac{\pi}{2n}} \left( \frac{2}{\pi}\, \right )^{\frac{1}{2}} \,dt =  \left( \frac{2}{\pi}\, \right )^{\frac{1}{2}}.
\]
Из полученных оценок $\, \|z^{n}\| \,$ следует, что в пространстве $\, BMOA \,\, \lim\limits_{n\rightarrow\infty} (\|z^{n}\|)^\frac{1}{n}=1$.

В пространствах $\mathcal{B}_{\alpha}$  при $n\geq 1$
\[
\|z^{n}\|= \sup\limits_{z\in D} |z^{n}|(1-|z|^{2})^{\alpha}\,= \left(\frac{n}{n+\alpha}\right)^{\frac{n}{2}} \,\left(\frac{2\alpha}{n+2\alpha}\right)^{\alpha},
\]
откуда следует, что  $\lim\limits_{n\rightarrow\infty} (\|z^{n}\|)^\frac{1}{n}=1$.

В пространстве $\mathcal{D}_{p}(\alpha)\,$  $ \|z^{n}\| = (\alpha_{n})^{\frac{1}{p}} $  и условиe теоремы 2.3 выполнено в случае, когда  $\lim \limits_{n\rightarrow\infty} (\alpha_{n})^\frac{1}{n}=1$.

Таким образом, условия теорем 2.1-2.3 выполнены во всех приведенных выше пространствах.  Для пространств $BMOA, \, \mathcal{B}_{\alpha}, \,\mathcal{D}_{p}(\alpha), \, \mathcal{A}_{p,q}^s (\mathbb{D})$  утверждения теорем 2.1-2.3 являются новыми.

\section{ Доказательства.}

Нам понадобятся следующие две леммы. Первая из них является естественным обобщением
неравенства Коши для коэффициентов ряда Тейлора.

\begin{lemma}
Пусть $ f\in X \quad$ и $ \quad
f(z)=\sum\limits_{k=0}^{\infty}c_{k}z^{k}\quad $в$\quad
\mathbb{D}.\quad$ Тогда
\[
|c_{n}|\,\|z^{n}\|\leq
E_{n}(f)\leq\|f\|.
\]

 \begin{proof}
Воспользовавшись формулой для коэффициентов ряда Тейлора запишем
равенство
\[
c_{n}z^{n}=\frac{1}{2\pi
i}\int\limits_{|\zeta|=1}\frac{f(z\zeta)-P_{n}(z\zeta)}{\zeta^{n+1}}d\zeta\quad,
\]
где $P_{n}$ - полином наилучшего приближения для  $f(z)$ степени не
выше $ (n-1)$. Отсюда $\quad|c_{n}|\|z^{n}\|\leq
E_{n}(f)\leq\|f(z)\|\quad$ ввиду свойств iii) и i)
 нормы в пространстве $X.$
 \end{proof}
\end{lemma}

\begin{lemma}
Пусть $ f\in X \quad$ и $\, \mu_{1}:=\liminf
\limits_{n\rightarrow\infty}\,(\|z^{n}\|)^{\frac{1}{n}}$,
$\mu_{2}:=\limsup
\limits_{n\rightarrow\infty}\,(\|z^{n}\|)^{\frac{1}{n}}$. Тогда
$ \mu_{1}\geq 1,  \,  \mu_{2} < \infty$.

 \begin{proof}
 Обозначим  $\, \beta_{n}=(\|z^{n}\|)^{\frac{1}{n}}$. Докажем, что  $\,  \mu_{2} < \infty$.
Предположим противное, тогда существует подпоследовательность  $\, \beta_{n_{k}}$  такая, что  $\lim
\limits_{k\rightarrow\infty}\, \beta_{n_{k}}=\infty $. Рассмотрим  функцию $ f_{0}$, определяемую
степенным рядом  $$ f_{0}(z)=\sum\limits_{k=0}^{\infty} \, (\beta_{n_{k}})^{\frac{-n_{k}}{2}}  z^{n_{k}}. \,$$
Она является целой и поэтому должна принадлежать $X$.  Но тогда согласно лемме 4.1  при каждом натуральном
$k \quad  (\beta_{n_{k}})^{\frac{-n_{k}}{2}} \| z^{n_{k}}\|\leq \|f_{0}\|<\infty \, $,  что невозможно.
Для доказательства того, что  $ \mu_{1}\geq 1,  \,$ допустим противное, т.е. что  $ \mu_{1}< 1$.
Зафиксируем $\varrho \in (\mu_{1}; \, 1)$  и рассмотрим функцию, определяемую рядом Тейлора

\begin{equation}
\label{trivial} \quad  f_{0}(z)=\sum\limits_{k=0}^{\infty} \varrho^{-n_{k}}z^{n_{k}},
\end{equation}
где ${n_{k}}$ последовательность, для которой $\liminf
\limits_{n\rightarrow\infty}\, \beta_{n}=\lim
\limits_{k\rightarrow\infty}\, \beta_{n_{k}}= \mu_{1}.$

Функция $f_{0}$  аналитична в круге $|z|<\rho,$ но неаналитична в  $\, \mathbb{D}$. Нетрудно проверить, что последовательность  частных сумм $S_{n, f_{0}}(z)$ ряда (10) фундаментальна в банаховом пространстве $X$ и, следовательно, сходится в нем к некоторой функции $f_{1} \in X$. Покажем, что тейлоровские коэффициенты функций $f_{0}$ и $f_{1}$  равны. Для произвольно фиксированного $k\in \mathbb{N}\bigcup \{ 0 \}, \, n>k $

\[
c_{k}(f_{1})=c_{k}(S_{n, f_{0}})+c_{k}(f_{1}-S_{n, f_{0}})=c_{k}(f_{0})+c_{k}(f_{1}-S_{n, f_{0}}).
\]
Устремив в этом равенстве $n \rightarrow \infty$  и учитывая лемму 2.1 получим  $c_{k}(f_{1})=c_{k}(f_{0})$. Таким образом, функция $f_{1} \in X$, но не является аналитической в $\mathbb{D}$, что противоречит характеристике пространства $X$. Следовательно, предположение о том, что  $ \mu_{1}< 1$  было неверным.

 \end{proof}
\end{lemma}

\begin{lemma}
Пусть $ f\in X $,
$\, K $- компакт, $\, K \subset \mathbb{D}$. Тогда  при $\, z \in K$

\[
|f(z)| \leq C \, \|f\|,
\]
где $C$ - постоянная, не зависящая  от $f$  и $z$.

 \begin{proof}
Положим  $\quad d:= \sup \{|z|: \quad z \in K\}, \,  d<1$. Запишем разложение функции $f$
в ряд Тейлора и оценим ее модуль, используя лемму 4.1 в предположении, что  $\, z \in K$.

\[
\quad
f(z)=\sum\limits_{k=0}^{\infty}c_{k}z^{k} ,
\]
\[
|f(z)|\leq \sum\limits_{k=0}^{\infty}|c_{k}|\,|z^{k}|\, \leq \,\|f(z)\| \sum\limits_{k=0}^{\infty}\frac{d^{k}}{\|z^{k}\|}\,\leq C \, \|f\|
\]
ввиду сходимости ряда, легко проверяемой с помощью признака Коши с учетом леммы 4.2.
\end{proof}
\end{lemma}

\begin{remark}
Из леммы 4.3 следует, что при выполнении ее условий значение функции в точке
$\, z \in K$ представляет собой ограниченный линейный фунционал и из сходимости
в пространстве $X$ следует равномерная сходимость на компактах  $ K \subset \mathbb{D}$.
Аналогичными рассуждениями можно показать, что при любом натуральном $n$  и $\, z \in K$
значение $f^{(n)}(z)$ является ограниченным линейным функционалом в $X$.
\end{remark}

Докажем теорему \ref{sol21}.
\begin{proof}
Достаточность.
Пусть $ f(z)=\sum\limits_{k=0}^{\propto}c_{k}z^{k}$
при $ z\in\mathbb{D}$.
\\ Согласно лемме 4.1 $ \quad |c_{n}|\,\|z^{n}\|\leq
E_{n}(f)$. Отсюда
\[
|c_{n}|\leq\frac{E_{n}(f)}{\|z^{n}\|}\quad
\mbox{и}\quad\lim\limits_{n\rightarrow
\infty}|c_{n}|^{\frac{1}{n}}\leq\lim\limits_{n\rightarrow
\infty}\left(\frac{E_{n}(f)}{\|z^{n}\|}\right)^{\frac{1}
{n}}=0\,,\quad\mbox{т.е. функция}\, f \mbox{- целая.}
\]

Необходимость. Для функции $f\in X$  положим  $f_\zeta (z):= f(z
\zeta)$. Поскольку $f$ - целая, то $f_{R}\in X$ при  любом $R>1$ и,
ввиду следствия 2 теоремы 2 \cite{DmChab} (см. также \cite{Dmz1},
теор. 2.3),
\[
E_{n}(f)\leq R^{-n}E_{n}(f_{R})\leq
R^{-n}\|f_{R}\|.
\]
Поэтому с учетом леммы 4.2
\[
0\leq\lim\limits_{n\rightarrow
\infty}\left(\frac{E_{n}(f)}{\|z^{n}\|}\right)^{\frac{1}{n}}
\leq \frac{1}{R}\,\limsup \limits_{n\rightarrow\infty}
\left(\frac{1}{\|z^{n}\|}\right)^{\frac{1}{n}}
\leq \frac{1}{R}.
\]
Ввиду леммы 4.2 и произвольности $R>1$  это означает, что
\[
\lim\limits_{n\rightarrow
\infty}\left(E_{n}(f) \right)^{\frac{1}{n}}=0.
\]
 \end{proof}

\begin{remark}
В работах   \cite{DmChab} и  \cite{Dmz1}  в условиях полученных там теорем пропущено условие iii),  без которого они, вообще говоря, неверны. В данной статье условие   iii)  предполагается выполненным и в этом случае используемое нами следствие 2 теоремы 2 работы \cite{DmChab} справедливо.
\end{remark}

Докажем теорему \ref{sol22}.
\begin{proof}
Достаточность. Из условия (8) следует, что выполнено условие теоремы
2.1 и, следовательно, $f$ - целая. Обозначим ее порядок $\alpha$ .
Тогда
 \begin{equation}
\alpha=\limsup \limits_{n\rightarrow\infty}\,\frac{n\ln
n}{-\ln|c_{n}|}\leq\limsup \limits_{n\rightarrow\infty}\,\frac{n\ln
n}{\ln\frac{\|z^{n}\|}{E_{n}(f)}}=\rho
\label{41}
 \end{equation}
ввиду леммы 4.1. Покажем теперь, что в условиях теоремы $\alpha>0$.
   Предположим противное, т.е. что
   $$
\limsup \limits_{n\rightarrow\infty}\,\frac{n\ln n}{-\ln|c_{n}|}=0.
   $$
   Тогда для произвольного положительного $\varepsilon\in (0,\,1)$ найдется $N_{\varepsilon}$  такое, что при
  $n>N_{\varepsilon}$  выполняется неравенство  $n\ln n<-\varepsilon \ln
  {|c_{n}|}$ и равносильное ему неравенство
  $$|c_{n}|<n^{\frac{-n}{\varepsilon}}.$$
  Пользуясь им оценим  $E_{n}(f)$  при $n> N_{\varepsilon}$.
Будем считать  $N_{\varepsilon}$ столь большим, что
$\|z^{n}\|\leq(\mu_{2}+\varepsilon)^{n}$ и
$\|z^{n}\|\geq(1-\varepsilon)^{n}$ при $n\geq
N_{\varepsilon}$. Тогда

\[
E_{n}(f)\leq\|\sum\limits_{k=n}^{\infty}c_{k}z^{k}\|\leq
\sum\limits_{k=n}^{\infty}k^{\frac{-k}{\varepsilon}}\,(\mu_{2}+\varepsilon)^{k}\leq\sum\limits_{k=n}^{\infty}
n^{\frac{-k}{\varepsilon}}\,\left(\mu_{2}+\varepsilon \right)^{k}=
\]
 \begin{equation}
=n^{\frac{-n}{\varepsilon}} (\mu_{2}+\varepsilon)^{n}(1-
\frac{\mu_{2}+\varepsilon}{n^{\frac{1}{\varepsilon}}})^{-1}
   \label{42}
 \end{equation}
при дополнительном условии $n>(\mu_{2}+\varepsilon)^{\varepsilon}.$
Из (10) получаем

\[
\frac{\|z^{n}\|}{E_{n}(f)}\geq
\left(\frac{1-\varepsilon}{\mu_{2}+\varepsilon}\right)^{n}\,n^{\frac{n}{\varepsilon}}\left(1-
\frac{\mu_{2}+\varepsilon}{n^{\frac{1}{\varepsilon}}}\right),
\]
\[
\ln
\left(\frac{\|z^{n}\|}{E_{n}(f)}\right)^{\frac{1}{n}}\geq
\ln
\frac{1-\varepsilon}{\mu_{2}+\varepsilon}+\,\frac{1}{\varepsilon}\ln
n +\frac{1}{n} \ln \left(1-
\frac{\mu_{2}+\varepsilon}{n^{\frac{1}{\varepsilon}}}\right).
\]
Отсюда
\[
\liminf \limits_{n\rightarrow\infty}\,
\frac{\ln\left(\frac{\|z^{n}\|}{E_{n}(f)}\right)^{\frac{1}{n}}}{\ln
n}\geq \frac{1}{\varepsilon},
\]
\[
\rho=\limsup \limits_{n\rightarrow\infty}\,\,\frac{n\ln n
}{\ln\frac{\|z^{n}\|}{E_{n}(f)}}\leq
\varepsilon,
\]
что противоречит условию теоремы.

 Выберем $\varepsilon\in(0;\frac{1}{2})\cap(0;\alpha)$. Из того, что
\[
\alpha=\limsup \limits_{n\rightarrow\infty}\,\frac{n\ln
n}{-\ln|c_{n}|}
\]
следует, что  существует $N_{\varepsilon}\in \mathbb{N}$, зависящее
только от  $\varepsilon$ и такое, что $|c_{n}|\leq
n^{-\frac{n}{\alpha+\varepsilon}}$  при всех $n\geq
N_{\varepsilon}$. Будем считать  $N_{\varepsilon}$ столь большим,
что $\|z^{n}\|\leq(\mu_{2}+\varepsilon)^{n}$ и
$\|z^{n}\|\geq(1-\varepsilon)^{n}$ при $n\geq
N_{\varepsilon}$. Тогда при $n> N_{\varepsilon}$
\[
E_{n}(f)\leq\|\sum\limits_{k=n}^{\propto}c_{k}z^{k}\|\leq
\sum\limits_{k=n}^{\propto}|c_{k}|\,\|z^{k}\|\leq\sum\limits_{k=n}^{\propto}
k^{\frac{-k}{\alpha+\varepsilon}}\,\|z^{k}\|\leq
\]
\begin{equation}
\label{L43}\leq\sum\limits_{k=n}^{\propto}n^{\frac{-k}{\alpha+\varepsilon}}(\mu_{2}+
\varepsilon)^{k}=\frac{(\mu_{2}+
\varepsilon)^{n}}{n^{\frac{n}{\alpha+\varepsilon}}}\cdot\left(1-\frac{\mu_{2}+
\varepsilon}{n^{\frac{1}{\alpha+\varepsilon}}}\right)^{-1}.
\end{equation}
Следовательно,
\[
\frac{\|z^{n}\|}{E_{n}(f)}\geq\frac{\|z^{n}\|}{(\mu_{2}+
\varepsilon)^{n}}\cdot
n^{\frac{n}{\alpha+\varepsilon}}\left(1-\frac{\mu_{2}+
\varepsilon}{n^{\frac{1}{\alpha+\varepsilon}}}\right),
\]
\begin{equation}
\label{L44}\alpha+\varepsilon\geq\frac{n\ln
n}{\ln\frac{\|z^{n}\|}{E_{n}(f)}}\cdot\left(1+\frac{\alpha+\varepsilon}
{n\ln n}\ln\left(1-\frac{\mu_{2}+
\varepsilon}{n^{\frac{1}{\alpha+\varepsilon}}}\right)+\frac{\alpha+\varepsilon}
{n\ln n}\ln\frac{\|z^{n}\|}{(\mu_{2}+
\varepsilon)^{n}}\right).
\end{equation}
Устремляя в (~\ref{L44}) $n\rightarrow \infty$, получим
$\alpha+\varepsilon\geq\rho$, что ввиду произвольности выбора
$\varepsilon>0$ означает, что $\alpha\geq\rho$. Таким образом,
$\alpha=\rho$ и достаточность доказана.

Необходимость. Пусть $f\in X$ - целая функция конечного
порядка $\rho$, т.е.
\begin{equation}
\label{L14}\limsup \limits_{n\rightarrow\infty}\,\frac{n\ln
n}{-\ln|c_{n}|}=\rho.
\end{equation}
Положим
\[
\alpha=\limsup \limits_{n\rightarrow\infty}\,\frac{n\ln
n}{\ln\frac{\|z^{n}\|}{E_{n}(f)}}
\]
($\alpha$ и $\rho$ в обозначениях по сравнению с доказательством
достаточности поменялись местами) и покажем, что $\alpha=\rho$. Из
леммы 4.1 аналогично (11) следует, что $\alpha\geq\rho$. Рассуждая
как при доказательстве достаточности мы можем утверждать, что для
произвольного $\varepsilon$,  $0<\varepsilon<1$ найдется
$N_{\varepsilon}$ такое, что $|c_{n}|\leq
n^{-\frac{n}{\rho+\varepsilon}}$ и $(1-\varepsilon)^{n}\leq
\|z^{n}\|\leq(\mu_{2}+\varepsilon)^{n}$ при $
n>N_{\varepsilon}$.
\\ Аналогично (~\ref{L43}) и (~\ref{L44}) (с заменой $\alpha$  на $\rho$ ) получим
\[
\rho+\varepsilon\geq\frac{n\ln
n}{\ln\frac{\|z^{n}\|}{E_{n}(f)}}\cdot\left(1+\frac{\rho+\varepsilon}
{n\ln n}\ln\left(1-\frac{\mu_{2}+
\varepsilon}{n^{\frac{1}{\rho+\varepsilon}}}\right)+\frac{\rho+\varepsilon}
{n\ln n}\ln\frac{\|z^{n}\|}{(\mu_{2}+
\varepsilon)^{n}}\right),
\]
откуда предельным переходом находим $\rho+\varepsilon\geq\alpha$ и,
следовательно, $\rho\geq\alpha$. Теорема доказана.
 \end{proof}

Докажем теорему \ref{sol23}.
\begin{proof}

 Достаточность. Пусть $f\in X$ и удовлетворяет
условию теоремы \ref{sol23} с некоторыми положительными $\rho$ и
$\sigma$. Тогда из (9) следует справедливость условия (\ref{L11})
теоремы \ref{sol22}, поэтому $f$ - целая и имеет порядок $\rho$. Пусть
тип $f$ равен $\alpha$. Докажем, что $\alpha=\sigma$. Из формулы для
определения типа целой функции

\begin{equation}
\label{L16} \alpha=\limsup
\limits_{n\rightarrow\infty}\,\frac{n}{e\rho}|c_{n}|
^{\frac{\rho}{n}}
\end{equation}
с учетом леммы 4.1 имеем $\alpha\leq\sigma$. Докажем обратное
неравенство. Из (\ref{L16}) следует, что для произвольного $\varepsilon>0$
существует $N_{\varepsilon}\in \mathbb{N}$ такое, что   при $n>
N_{\varepsilon}$
\begin{equation}
\label{L17}|c_{n}| <\left(\frac{\rho
e(\alpha+\varepsilon)}{n}\right)^{\frac{n}{\rho}}.
\end{equation}
С учетом (17) подобно (13), (14) находим
\[
E_{n}(f)\leq\sum\limits_{k=n}^{\infty}\left(
\frac{\rho
e(\alpha+\varepsilon)}{k}\right)^{\frac{k}{\rho}}\|z^{k}\|\leq
\left( \frac{\rho
e(\alpha+\varepsilon)}{n}\right)^{\frac{n}{\rho}}(\mu+\varepsilon)^{n}\times
\]
\begin{equation}
\label{L18}\times\left(1-\frac{C}{n^{\frac{1}{\rho}}}\right)^{-1}
\quad,
\end{equation}
где $C=(\mu+\varepsilon)(\rho
e(\alpha+\varepsilon))^{\frac{1}{\rho}}$. Из (\ref{L18}) получаем
\begin{equation}
\label{L19}\alpha+\varepsilon\geq\frac{n}{e\rho}\left(\frac{E_{n}(f)}
{\|z^{n}\|}\right)
^{\frac{\rho}{n}}\frac{\|z^{n\|^{\frac{\rho}{n}}}}{(\mu+\varepsilon)^{\rho}}\left(1-\frac{c}
{n^{\frac{1}{\rho}}}\right)^{\frac{\rho}{n}}.
\end{equation}
Вычислив верхний предел в (\ref{L19}), имеем
\begin{equation}
\label{L20}\alpha+\varepsilon\geq\sigma\left(\frac{\mu}{\mu+\varepsilon}\right)^{\rho},
\end{equation}
откуда, устремляя $\varepsilon$ к нулю, находим $\alpha\geq\sigma$ ,
что завершает доказательство достаточности.

Необходимость. Пусть $f\in X$ - целая, конечного порядка и
нормального типа. Обозначим ее порядок $\rho$ (он удовлетворяет
ввиду теоремы \ref{sol22} равенству (8)) и обозначим $\alpha$ ее
тип. Покажем, что $\alpha=\sigma$. Ввиду (\ref{L16}) и леммы 4.1
 $\alpha\leq\sigma$. Далее для доказательства
неравенства $\alpha\geq\sigma$ нужно полностью повторить
соответствующие рассуждения из доказательства достаточности.
 \end{proof}

В теоремах 2.1 - 2.3 предполагается, что функция $f\in X$ является
аналитической в единичном круге. Это требование можно снять, несколько изменив формулировки теорем.
Если считать, что $X$ - банахово пространство функций, определенных в круге $\mathbb{D}$,  норма в котором удовлетворяет условиям  i), ii) и  iii) (в отличие от предшествующего рассмотрения мы не требуем от  $f\in X$ аналитичности в $\mathbb{D}$), то будут справедливы следующие модификации теорем 2.1 - 2.3.

\begin{theorem}\label{sol45}
    Пусть  функция $f\in X$ и
$\liminf \limits_{n\rightarrow\infty}\,(\|z^{n}\|)^{\frac{1}{n}}=\mu_{1}>0$.
Если $f$ - целая, то
 \begin{equation}
 \label{L45}
\lim\limits_{n\rightarrow
\infty}\,\left\{\frac{E_{n}(f)}{\|z^{n}\|}\right\}^{\frac{1}{n}}=0.
 \end{equation}
Обратно, если выполнено условие (\ref{L45}), то последовательность $\{p^{*}_{n}\}$
полиномов наилучшего приближения функции $f$ по норме пространства  $X$
равномерно сходится в каждом круге $|z| < r, \, r\in (0; \mu_{1})$
к некоторой целой функции.
           \label{45}
      \end{theorem}

\begin{proof}
Первое утверждение доказано в теореме 2.1. Докажем второе.
Пусть для $f\in X$ выполнено условие (\ref{L45}). Из него следует, что
$\, E_{n}(f)\rightarrow 0$ при  $\, n \rightarrow \infty.$
Кроме того, для любых натуральных $m$ и $n, \, m \geq n$ и
последовательности $\{P^{*}_{n}(z)\}$  полиномов наилучшего приближения
$\|P^{*}_{n}(z)-P^{*}_{m}(z)\| \leq 2E_{n}(f)$.
Из леммы 4.3 следует, что последовательность полиномов наилучшего приближения
фундаментальна по  $sup$ - норме в  каждом круге $|z|\leq r$  при $r \in
(0; \mu_{1})$  и справедлива оценка $|P^{*}_{n}(z)-P^{*}_{m}(z)| \leq 2CE_{n}(f)$
с постоянной $C$, зависящей только от $r$ и $\mu_{1}$. Поэтому  последовательность
 $\{P^{*}_{n}(z)\}$ сходится равномерно на компактах в  круге  $ |z|< \mu_{1}$ к
 некоторой функции  $g(z)$, аналитической в круге $ |z|< \mu_{1}$, причем
 $|P^{*}_{n}(z)-g(z)| \leq 2CE_{n}(f)$  при $|z|\leq r$. Из неравенств Коши для
 коэффициентов $\gamma_n$ ряда Тейлора функции $g(z)$ следует оценка
 $|\gamma_{n}| \leq 2CE_{n}r^{-n}(f)$. Для  доказательства аналитичности функции $g(z)$
 во всей плоскости достаточно воспользоваться формулой Коши-Адамара для радиуса
 сходимости степенного ряда.

 \end{proof}

Формулировки теорем 2.2 и 2.3 изменяются следующим образом.

\begin{theorem}\label{sol46}
Пусть $\liminf
\limits_{n\rightarrow\infty}\,(\|z^{n}\|)^{\frac{1}{n}}=\mu_{1}>0
\,,\quad\limsup
\limits_{n\rightarrow\infty}\,(\|z^{n}\|)^{\frac{1}{n}}=\mu_{2}<\infty$,  $f\in X$.

\ Если функция $f$ - целая конечного порядка
$\rho\in(0;\infty)$, то
\begin{equation}
\label{L46}\limsup \limits_{n\rightarrow\infty}\,\frac{n\ln
n}{\ln\frac{\|z^{n}\|}{E_{n}(f)}}=\rho .
\end{equation}

Обратно, если выполнено условие (\ref{L46}), то последовательность $\{p^{*}_{n}\}$
полиномов наилучшего приближения функции $f$ по норме пространства  $X$
равномерно сходится в каждом круге $|z| < r, \, r\in (0; \mu_{1})$
к некоторой целой функции конечного порядка $\rho\in(0;\infty)$, определяемого равенством (\ref{L46}).
 \label{46}
    \end{theorem}

\begin{theorem}\label{sol47}
Пусть существует конечный предел $\lim\limits_{n\rightarrow \infty}
(\|z^{n}\|)^\frac{1}{n}=\mu>0 $,  $f\in X$. Если $f$ -  целая функция конечного порядка $\rho\in(0;\infty)$ и
нормального типа $\sigma\in(0;\infty)$, то
\begin{equation}
\label{L47}\limsup\limits_{n\rightarrow\infty}
\,\frac{n}{e\rho}\left(\frac{E_{n}(f)}{\|z^{n}\|}\right)
^{\frac{\rho}{n}}=\sigma .
\end{equation}
Обратно, если выполнено условие (\ref{L47}), то последовательность $\{p^{*}_{n}\}$
полиномов наилучшего приближения функции $f$ по норме пространства  $X$
равномерно сходится в каждом круге $|z| < r, \, r\in (0; \mu_{1})$
к некоторой целой функции конечного порядка $\rho\in(0;\infty)$
нормального типа $\sigma\in(0;\infty)$ ( $\sigma, \rho$ из соотношения (\ref{L47})).

\label{47}
 \end{theorem}

\begin{remark}
Условия $\mu_{1}>0$,  $\mu_{2}<\infty$ в теоремах 4.1-4.2 достаточно
естественны (легко привести примеры, показывающие, что без этих
условий соответствующие теоремы неверны).  Нельзя ли заменить ими
условие существования конечного предела $\lim\limits_{n\rightarrow
\infty} (\|z^{n}\|)^\frac{1}{n}=\mu>0 $ в теореме 4.3? Этот вопрос
остается открытым.
\end{remark}

\section{ О возможности аппроксимации функций $f\in X$ полиномами из ${\mathcal{P}}_{n}[\mathbb{Z}]$ }

Вопрос о возможности равномерного приближения  функций, непрерывных
на компактах в $\mathbb{C}$ и аналитических во внутренних точках
с любой точностью многочленами с целыми  коэффициентами
 хорошо изучен. Историю вопроса и полученные результаты можно найти
 в статье С.Я. Альпера \cite{Alp} (см. также статью обзорного характера
 \cite{Trig}). В частности, для того чтобы отличную от полинома аналитическую в круге
 с центром  $z=0$ и непрерывную на его замыкании функцию можно было
 равномерно аппроксимировать с любой точностью, необходимо, чтобы радиус круга
 был меньше единицы.
 В других пространствах вопрос о возможности аппроксимации аналитических
 функций многочленами из  ${\mathcal{P}}_{n}[\mathbb{Z}]$ менее изучен.
 Мы рассмотрим ситуацию с возможностью аппроксимации
 аналитических функций многочленами с целыми  коэффициентами в изучаемых
 в данной статье пространствах $X$.  Будем далее считать, как и всюду в статье,
 кроме теорем  4.1 - 4.4, что в $X$ содержатся только функции, аналитические в  $\mathbb{D}$.

\begin{theorem}\label{sol51}

Пусть аналитическая в единичном круге функция $f\in X$ и
$\inf \limits_{n \in \mathbb{N}}\, \|z^{n}\|>0$. Если существует
последовательность $p_{n} \in {\mathcal{P}}_{n}[\mathbb{Z}]$
такая, что
$$
\lim\limits_{n\rightarrow
\infty}\,\|f-p_{n}\|=0,
$$
то $f$ является полиномом с целыми коэффициентами.
   \end{theorem}
Подобная теорема  для пространств $B$ и $H_p$ известна \cite{Igor}.
     Докажем теорему \ref{sol51}.
\begin{proof}
$$
\|p_{n+1} - p_{n}\| \leq \|p_{n+1} - f\| + \|f - p_{n}\| \rightarrow 0
$$
при $n \rightarrow \infty$. С другой стороны, если $p_{n+1} \neq p_{n}$
тождественно, то разность $p_{n+1} - p_{n}$ есть ненулевой многочлен с
целыми коэффициентами и вследствие леммы 4.1

$$
\|p_{n+1} - p_{n}\| \geq \inf \limits_{n \in \mathbb{N}}\, \|z^{n}\| ,
$$
т.е. $\|p_{n+1} - p_{n}\|$ не стремится к нулю. Следовательно, при
всех $n \quad p_{n+1} \equiv p_{n} \equiv f$.

 \end{proof}

Следующее утверждение хорошо известно в случае аппроксимации в
пространстве $B$  и при равномерной сходимости на компактах в
единичном круге (см. например, \cite{Igor}).

\begin{theorem}\label{sol52}

Пусть аналитическая в единичном круге функция $f \in X$. Если существует
последовательность $p_{n} \in {\mathcal{P}}_{n}[\mathbb{Z}]$
такая, что
$$
\lim\limits_{n\rightarrow
\infty}\,\|f-p_{n}\|=0,
$$
то все ее тейлоровские коэффициенты  $c_{k}:= \frac{f^{k}(0)}{k!}$ - целые числа.

      \end{theorem}

\begin{proof}
Многочлен  $p_{n}$ можно записать в форме

$$
p_{n}(z) = \sum\limits_{k=0}^{n}a_{k,n}z^{k},
$$
где коэффициенты   $a_{k,n}$  целые. Ввиду леммы 4.1

$$
|c_{k} - a_{k,n}| \,  \|z^{k}\| \leq \|f(z) - p_{n}(z)\| \rightarrow 0 , \quad n \rightarrow \infty .
$$
Следовательно,
$$
c_{k}= \lim\limits_{n\rightarrow
\infty}\,a_{k,n}\in \mathbb{Z}.
$$

 \end{proof}

      \begin{corollary}
Из теоремы 5.2 следует, что при выполнении ее условий радиус сходимости
ряда Тейлора по степеням $\, z$ для функции, отличной от полинома,
не больше единицы.
\end{corollary}

\begin{theorem}\label{sol53}

Пусть $ f \in X \quad$,  $ \quad
f(z)=\sum\limits_{k=0}^{\infty} c_{k}z^{k}\quad $в$\quad
\mathbb{D}\quad$  и
$$
\lim\limits_{n\rightarrow
\infty} \, \|f(z)-S_{n}(f,z)\|=0,
$$
где $ S_{n}(f,z)$ - частная сумма ряда Тейлора функции $f$ порядка $n$.

Для существования последовательности $p_{n} \in {\mathcal{P}}_{n}[\mathbb{Z}]$
такой, что
$$
\lim\limits_{n\rightarrow
\infty}\,\|f(z)-p_{n}(z)\|=0,
$$
необходимо и достаточно, чтобы все коэффициенты $c_k$ были целыми.

      \end{theorem}

\begin{proof}
Необходимость следует из теоремы 5.2, достаточность - из того, что
$ S_{n}(f,z)$ есть многочлен с целыми коэффициентами.

 \end{proof}

Следующая теорема проясняет различие между пространствами $H_p$ и
$H'_{p}, \quad p\in (1, \infty)$,  обнаруженное И. Прицкером в \cite{Igor}
( в пространствах $H'_{p}, \quad p\in (1, \infty)$ нетривиальные функции
с целыми тейлоровскими коэффициентами можно аппроксимировать полиномами из
$ {\mathcal{P}}_{n}[\mathbb{Z}]$,  а в  пространстве $H_p$ - нет).

\begin{theorem}\label{sol54}

Пусть для некоторой аналитической в единичном круге функции $ f \in X, \quad$
отличной от полинома,
существует последовательность полиномов  $p_{n} \in {\mathcal{P}}_{n}[\mathbb{Z}]$
такая, что

$$
\lim \limits_{n \rightarrow
\infty} \,\|f(z)-p_{n}(z)\|=0.
$$
Тогда
$\liminf \limits_{n \rightarrow \infty} \, \|z^{n}\|=0.$

      \end{theorem}

 \begin{proof}
Сопоставим каждому полиному $p_{n}$ из условия теоремы функцию $z^{m_{n}}$, выбирая число $m_{n}$ так, чтобы  $m_{n} > deg \, p_{n}$, коэффициент $c_{m_{n}}$ в тейлоровском разложении функции был отличен от нуля и последовательность $m_{n}$ строго возрастала (понятно, что такую последовательность можно выбрать многими способами; нам подходит любая). Поскольку ввиду теоремы 5.2 коэффициенты  $c_{m_{n}}$  целые, то
$$
\| z^{m_{n}}\| \leq | c_{m_{n}}| \, \| z^{m_{n}}\| \leq \|f(z)-p_{n}(z)\| \rightarrow 0
$$
при  $n \rightarrow \infty$.

 \end{proof}

Необходимое условие   $\liminf \limits_{n \rightarrow \infty} \, \|z^{n}\|=0$  теоремы  \ref{sol54} является также в некотором смысле достаточным, что показывает следующая теорема.

\begin{theorem}\label{sol55}

Для того, чтобы в пространстве $ X$ существовала функция  $f$,  отличная от полинома
и сколь угодно хорошо аппроксимируемая  последовательностью полиномов
$p_{n} \in {\mathcal{P}}_{n}[\mathbb{Z}]$,  необходимо и достаточно, чтобы  в
пространстве $ X$  выполнялось условие   $\liminf \limits_{n \rightarrow \infty} \, \|z^{n}\|=0.$

      \end{theorem}

\begin{proof}
Необходимость следует из предыдущей теоремы. Докажем достаточность.  Из условия  $\liminf \limits_{n \rightarrow \infty} \, \|z^{n}\|=0$ следует,  что существует возрастающая последовательность натуральных чисел $n_{k}$  такая, что $\| z^{n_{k}}\| \leq 2^{-k}$  при всех  $k = 1, 2, 3,... \,$. Положим

 $$
 f(z)=\sum\limits_{k=0}^{\infty} z^{n_{k}}\quad.
$$
 Нетрудно видеть, что $f$   аналитична в  единичном круге,  $ f \in X \,$  и $f$  сколь угодно хорошо аппроксимируема полиномами с целыми коэффициентами  (частными суммами  $ S_{n}(f,z)$ ).

 \end{proof}

\begin{remark}
Заметим, что полученные в статье утверждения остаются верными и в случае, если условие
 iii) заменить более слабым условием
\begin{equation}
iii^{\prime}) \quad \, \|  \frac {1}{2\pi}\int\limits_{0}^{2\pi}f(ze^{it})g(t)
\,dt \| \leq \, C \, \frac {1}{2\pi}\int\limits_{0}^{2\pi}|g(t)|
\,dt \,\, \|f(\cdot)\|,
\end{equation}
где постоянная $C$ не зависит от  $f\in X$  и $g\in L_{[0; \, 2\pi]}$.

\end{remark}

Сведения об авторе:

Двейрин Михаил Захарович

 Донецкий национальный университет,

Кафедра математического анализа и теории функций

г. Донецк, 83055 Ул. Университетская, 24

Е-mail: matem47@mail.ru

тел. +38(062)-2972953, +380667075599


\begin{thebibliography}{1}


\bibitem {HL2}
 Hardy G.H., Littlewood J.E. \emph{Some properties of fractional integrals II}
Math. Z., 1931, 34, N 3, 403--439.

\bibitem {DRS}
 Duren P.L., Romberg B.W., Shields A.L. \emph{Linear functionals in
$H_p$  spaces with  $0<p<1$} // J. reine und angew. Math. - 1969. - 238,
s. 4--60.

\bibitem{Gvar}
Гварадзе М. И.   \emph{Об одном классе пространств аналитических
функций. } Мат. заметки // - 1977. 21, N 2, с. 141--150.

\bibitem{Swed}
Шведенко С. В.  \emph{Классы Харди и связанные с ними пространства
аналитических функций в единичном круге и шаре.} // Итоги науки и
техники. Сер. мат. анализ, Москва, ВИНИТИ - 1985. 23 - с. 3--124.

\bibitem {Zhu}
 K. Zhu \emph{Bloch type spaces of analytic functions}// Rocky mountain J. Math. - 1993. - 23, N 3, p. 1143--1177.

\bibitem{Dink}
Дынькин Е. М.  \emph{Конструктивная характеристика классов С.Л.
Соболева и О.В. Бесова.}// Труды мат. ин-та АН СССР - 1981. - 155,
с. 41--76.

\bibitem{Rd1}
Reddy A.R. \emph{A Constribution to best approximation in the
$L^{2}$ norm.} // J. Approxim. Theory. 1974. - 11, N 11, p.
110--117.

\bibitem{Ibr1}
 Ибрагимов И.И., Шихалиев Н.И.  \emph{О наилучшем полиномиальном
приближении в одном пространстве аналитических функций.}// ДАН СССР. - 1976. - 227, N 2, с. 280--283.

\bibitem{Ibr2}
 Ибрагимов И.И., Шихалиев Н.И.  \emph{О наилучшем
 приближении в среднем  аналитических функций в пространстве $A_{p}(|z|<1).$
} // Спец. вопросы теории функций.-Баку: ЭЛМ -1977. N 1, с. 84--96.

\bibitem{Wac1}
 Вакарчук С.Б.  \emph{О наилучшем полиномиальном
приближении аналитических функций в  пространстве
$\mathcal{B}_{p,\,q,\,\lambda}.$} // ДАН УССР, сер. физ.-мат. и техн.
науки. - 1989. - N 8, с. 6--9.

\bibitem{Wac2}
 Вакарчук С.Б.  \emph{О наилучшем полиномиальном
приближении аналитических в единичном круге функций.}// УМЖ. -1990. - 42,  N 6, с. 838--843.

\bibitem{Mamad}
 Мамадов Р.  \emph{Некоторые вопросы приближения целыми функциями.}
 Автореферат диссертации на соискание ученой степени
кандидата физико-математических наук,спец. 01.01.01 - математический анализ. - Душанбе, 2009. -
14с.

\bibitem{DmChab}
 Двейрин М.З., Чебаненко И.В.  \emph{О полиномиальной аппроксимации в банаховых пространствах
аналитических функций.}  // Теория отображений и приближение
функций. -Киев, Наукова думка. - 1983, с. 63--73.

\bibitem{Dmz1}
 Двейрин М.З. \emph{Неравенство Адамара и наилучшее приближение
функций, аналитических в единичном круге} // Укр. матем. вестник. -
2006. - 3, N 3, с. 315--330.


\bibitem{Alp}
Альпер С.Я.   \emph{О приближении функций многочленами с целыми коэффициентами на замкнутых множествах. }// Изв. АН СССР, сер. матем. - 1964. - 28, N 5 - с. 1173--1186.

\bibitem{Trig}
Тригуб Р.М. \emph{Приближение функций с диофантовыми условиями многочленами с целыми коэффициентами.}  // Метрические вопросы теории функций и отображений. -Киев, Наукова думка. - 1971, с. 267--333.

\bibitem{Igor}
Igor E. Pritsker \emph{An areal analog of Mahler's measure.}//
Illinois J. Math. - 2008. - 52, N 2 - p. 347-363.


\end{thebibliography}
\end{document}